\documentclass[12pt]{amsart} 

\newtheorem{theorem}{Theorem}[section]
\newtheorem{proposition}[theorem]{Proposition}
\newtheorem{definition}[theorem]{Definition}
\newtheorem{lemma}[theorem]{Lemma}

\newtheorem{remark}[theorem]{Remark}

\newcommand{\Hn}{\hat{0}}
\newcommand{\He}{\hat{1}}
\newcommand{\Z}{\mathbb{Z}}
\newcommand{\N}{\mathbb{N}}
\newcommand{\oI}{\overline{I}}
\newcommand{\oP}{\overline{P}}
\newcommand{\Br}{\mathrm{Br}}
\newcommand{\In}{\mathrm{Invol}}
\newcommand{\tW}{\widetilde{W}}
\newcommand{\tS}{\widetilde{S}}
\newcommand{\lW}{W^-}
\newcommand{\ts}{\widetilde{s}}
\newcommand{\twist}{\mathfrak{I}}
\newcommand{\id}{\mathrm{id}}
\newcommand{\inv}{\mathrm{inv}}

\begin{document}

\title[Fixed points]{Fixed points of involutive automorphisms of the
  Bruhat order}  
\thanks{Partially supported by the European Commission's IHRP
  Programme, grant HPRN-CT-2001-00272, ``Algebraic Combinatorics in
  Europe''.}

\author{Axel Hultman}
\address {Fachbereich Mathematik und Informatik, Philipps-Universit\"at Marburg, Hans Meerwein Strasse, Lahnberge, D-35032 Marburg, Germany}
\email{axel@mathematik.uni-marburg.de} 

\begin{abstract}
Applying a classical theorem of Smith, we show that the poset property of
being Gorenstein$^\ast$ over $\Z_2$ is inherited by the subposet of
fixed points under an involutive poset automorphism. As an application,
we prove that every interval in the Bruhat order on (twisted) involutions in an
arbitrary Coxeter group has this property, and we find the rank
function. This implies results conjectured by F.\ Incitti. We also
show that the Bruhat order on the fixed points of
an involutive automorphism induced by a Coxeter graph automorphism is
isomorphic to the Bruhat order on the fixed subgroup viewed as a
Coxeter group in its own right.
\end{abstract} 

\maketitle

\section{Introduction}
In \cite{RS,RS2}, Richardson and Springer initiated the study of Bruhat
decompositions of certain symmetric varieties. They carried out the
following construction. Consider a connected, reductive linear
algebraic group $G$ over an algebraically closed field $F$ with
${\rm char}(F) \neq 2$. Let $B \subseteq G$ be a Borel subgroup and $T
\subseteq B$ a maximal torus. Given a $G$-automorphism $\theta$ of
order $2$ preserving $T$ and $B$, let $K$ be the fixed point
group. Define the {\em symmetric 
  variety} $X = G/K$. Now, $B$ acts by left translations on $X$, giving rise to a finite number of orbits. We may order these orbits by containment of their Zariski closures. The following special case is worth mentioning: $G$ is a symmetric variety for $G \times G$, and the orbits under the $B\times B$-action coincide with the $B$-orbits of the flag variety $G/B$. In this case, the order obtained is the Bruhat order on the corresponding Weyl group. 

The way in which Richardson and Springer studied this order was by
means of an order-preserving map to the subposet of twisted
involutions in the Bruhat order on $W$. When $\theta$ acts trivially
on $T$, this is just the Bruhat order on the involutions of $W$. The
latter poset has been studied by Incitti \cite{incitti, incitti2, incitti3} who
showed that it is EL-shellable (hence Cohen-Macaulay) and Eulerian when $W$ is 
a classical Weyl group. In \cite{incitti3}, these properties were
conjectured to hold for arbitrary Coxeter groups. (In infinite groups,
this should be interpreted as these properties holding for every interval.) He
also predicted an interpretation for the rank function. 

In this paper, we prove that every interval in the Bruhat order on the twisted
involutions of an arbitrary Coxeter group (with respect to an
involutive group automorphism which preserves the Coxeter generator set) is
Gorenstein$^\ast$ over $\Z_2 = \Z/2\Z$. Moreover, we find the rank
function. This implies ``most'' of
Incitti's conjecture, namely Eulerianity and Cohen-Macaulayness
over $\Z_2$, as well as the assertion about the rank function.

We also study the subposet of Bruhat order induced by the
fixed subgroup (actually a Coxeter group) of an involutive group automorphism
preserving the set of Coxeter generators. The subposet turns out to be
isomorphic to this subgroup's own Bruhat order.

Both results rely heavily on a general theorem stating that the
Goren\-stein$^\ast$ property over $\Z_2$ is inherited by induced
subposets of fixed points under involutive poset automorphisms. This
is proved using one of Smith's classical results on group actions on spheres.

The remainder of the paper has the following structure. We review some
preliminaries on combinatorial topology and Coxeter groups in Section
\ref{se:prel}. In Section \ref{se:main}, we recall the classical
result of Smith referred to above. We use it to prove the result which
forms the technical backbone of the paper, namely that the
Gorenstein$^\ast$ property over $\Z_2$ is inherited by subposets of
fixed points under involutive automorphisms. The rest of the
paper is devoted to applications of this theorem to Bruhat
orders. In Section \ref{se:twist}, we study Bruhat orders on twisted
involutions, and in Section \ref{se:graphauto}, we focus on induced
Bruhat orders on fixed point subgroups of automorphisms induced
by Coxeter graph automorphisms. In the latter section, related results
for the weak order are also proved. 

We are grateful to an anonymous referee who pointed out to us a way to generalize the main result of Section \ref{se:graphauto} using methods different from ours. The referee's argument is stated in Appendix \ref{se:referee}. 

\section{Preliminaries}\label{se:prel}

\subsection{Posets and simplicial complexes}
Let $\Delta$ be a finite simplicial complex. Its {\em dimension} is
the maximum dimension of a {\em facet}, i.e.\ inclusion-maximal face.
We say that $\Delta$ is {\em pure} if its facets are equidimensional. 

Suppose $\Delta$ is pure of dimension $n$. Let $\mathcal{F}$ be its set
of facets. Then $\Delta$ is {\em strongly connected} if for any pair
$F,G \in \mathcal{F}$, there exists a sequence of facets $F=F_0, F_1,
\dots, F_t=G$ such that $F_{i-1} \cap F_i$ has dimension $n-1$ for all
$i \in [t] = \{1, \dots, t\}$. We call $\Delta$ {\em thin} if every
$(n-1)$-dimensional face is contained in exactly two facets.
\begin{definition}
A finite simplicial complex is a {\em pseudomanifold} if it is pure,
thin and strongly connected.
\end{definition}
It is easy to see that if $\Delta$ is an $n$-dimensional
pseudomanifold, then $\widetilde{H}_n(\Delta;\Z_2) \cong \Z_2$, where
$\widetilde{H}_\ast$ denotes reduced homology. This is since, over
$\Z_2$, the (homology class of the) sum of all facets is the unique
nontrivial element in the top-dimensional reduced simplicial homology.

A poset $P$ is {\em bounded} if it has unique top and bottom elements,
denoted $\He$ and $\Hn$, respectively. If $P$ is bounded, then its
{\em proper part} is $\oP = P\setminus\{\Hn,\He\}$. 

\begin{definition}
A poset $P$ is {\em Eulerian} if it is bounded, graded and finite, and
  its M\"obius function satisfies $\mu(p,q) =
  (-1)^{\rho(q)-\rho(p)}$ for all $p \leq 
q \in P$, where $\rho$ is the rank function of $P$.
\end{definition}

To any poset $P$, we may associate the {\em order complex} $\Delta(P)$. This
is the simplicial complex whose faces are the chains in $P$. Thus, we
can assign topological properties to $P$. If $P$ is bounded, however,
$\Delta(P)$ is not very exciting; the extremal elements serve as cone
points in the complex. Thus, when we speak of topological properties
of a bounded poset $P$, we have the corresponding properties of
$\Delta(\oP)$ in mind. 

We say that $P$ has the {\em diamond property} if every
interval of length $2$ in $P$ is isomorphic to the ``diamond-shaped''
four-element poset consisting of two incomparable elements, together
with a top and a bottom element. Note that a graded poset is thin iff it is
has the diamond property.

In the definitions that follow, $k$ is any abelian group. We will
primarily be interested in the case $k = \Z_2$.
\begin{definition}
A poset $P$ is {\em
  Cohen-Macaulay over $k$} if it is bounded, graded and finite, and
  every interval $[p,q] \subseteq P$ satisfies
  $\widetilde{H}_i([p,q];k) = 0$ for all $0 \leq i <
  \rho(q)-\rho(p)-2$, and $\widetilde{H}_{\rho(q)-\rho(p)-2}([p,q];k)
  \cong k^m$ for some $m$, where $\rho$ is the rank function of $P$.   
\end{definition}

In other words, for $P$ to be Cohen-Macaulay, the homology of every
interval must be the homology of a (possibly empty) wedge of top-dimensional
spheres. We may restrict this further to require the number of spheres
to be exactly one, yielding the following definition:

\begin{definition}
A poset $P$ is {\em Gorenstein$^\ast$ over $k$} if it is bounded,
  graded and finite, and every interval is a top-dimensional homology
  sphere over $k$.
\end{definition}

Using the correspondence between the M\"obius function and the Euler
characteristic (Ph.\ Hall's Theorem), one obtains the following
alternative definition: 
\begin{proposition}
A poset is Gorenstein$^\ast$ over $k$ iff it is Cohen-Macau\-lay over
$k$ and Eulerian. 
\end{proposition}

Sometimes we refer to a poset as being Cohen-Macaulay
(Goren\-stein$^\ast$) without declaring over which group. We then have
the group $\Z$ in mind. By the Universal Coefficient Theorem, this
implies the Cohen-Macaulay (Gorenstein$^\ast$) property over any
abelian group. 

\subsection{Some facts about Coxeter groups} 
Here, we collect some facts that we need about Coxeter groups and their partial
orders. We refer to Humphreys \cite{humphreys} for a thorough
background to these matters. 

Throughout the rest of the paper, $(W,S)$ will be a Coxeter system
with $|S| < \infty$ and length function $\ell:W\to \N$. We will always
assume that symbols of the form $s_i$ are elements in $S$. If $w =
s_1\dots s_k$ and $\ell(w) = k$, then $s_1 \dots s_k$ is called a {\em
  reduced expression} for $w$. Two important properties (that in fact
are equivalent and characterize Coxeter systems) are the following: 

\begin{proposition}[Deletion Property] Suppose $s_1 \dots s_k$ is a
  non-redu\-ced expression for $w$. Then there exist indices $1\leq i <
  j \leq k$ such that $s_1\dots \widehat{s_i} \dots \widehat{s_j}
  \dots s_k = w$, where the hats denote omission.
\end{proposition}

\begin{proposition}[Exchange Property] Let $s_1\dots s_k$ be any
  expression for $w$. If $\ell(w) > \ell(ws)$ for some $s \in S$, then
  $ws = s_1\dots \widehat{s_i} \dots s_k$ for some $i \in [k]$.
\end{proposition}

\begin{remark}
Let $T$ denote the set of reflections in $W$. Replacing the hypothesis
$s \in S$ by the weaker $s \in T$ in the 
statement of the Exchange Property yields
another true statement known as the {\em Strong Exchange Property}.
\end{remark}

If we are interested in the set of all reflections rather than the set
of simple reflections $S$, we may define the {\em absolute length}
$\ell^\prime:W \to \N$. Then, $\ell^\prime(w)$ is the smallest $k$
such that $w$ is a product of $k$ reflections. Clearly,
$\ell^\prime(w) \leq \ell(w)$ for all $w \in W$.

We now define the two most important ways to partially order $W$. 

\begin{definition} The (right) {\em weak order} on $W$ is defined by $u \leq
  v$ iff $v = uw$ and $\ell(v) = \ell(u) + \ell(w)$ for some $w \in W$.
\end{definition}

Clearly, the length function $\ell$ serves as rank function of the
weak order. 

\begin{definition}
The {\em Bruhat order} on $W$ is defined by $v \leq
w$ iff some (equivalently, every) reduced expression $s_1\dots s_k$ for
$w$ contains a subexpression $s_{i_1}\dots s_{i_j}$, $1 \leq i_1 <
\dots < i_j \leq k$, which is a reduced expression for $v$. We denote this
poset by $\Br(W)$.
\end{definition}

It is obvious that the Bruhat order contains the weak order as
relations. Although 
not immediate from the definition, $\Br(W)$, too, is graded with
rank function $\ell$. Clearly, every interval in $\Br(W)$
is finite, even if $W$ is infinite. Moreover, the intervals have a
nice topological structure: 

\begin{theorem}[Bj\"orner and Wachs \cite{bjornerwachs}] \label{th:BW}
Given any Coxeter group $W$, every interval in $\Br(W)$ is homeomorphic
to a sphere of top dimension.  
\end{theorem}

Any (labelled) graph automorphism of the Coxeter graph of $W$
of course induces an automorphism of $\Br(W)$. A (slightly) less trivial
automorphism of the latter is given by the inversion map $w \mapsto
w^{-1}$. Since it leaves all $s \in S$ fixed, but not all $w \in W$
(in general), it cannot be induced by a graph automorphism.

The dihedral groups are easy to deal with separately, but they do not fit into
the following picture:
\begin{theorem}[van den Hombergh \cite{hombergh}, Waterhouse
    \cite{waterhouse}]\label{th:hombergh} 
If $W$ is irreducible and $|S| > 2$, the automorphism group of $\Br(W)$ is generated by $w \mapsto w^{-1}$ and the automorphisms induced by Coxeter graph automorphisms. 
\end{theorem}

\section{Fixed points of poset automorphisms}\label{se:main}
Consider an involutive automorphism (i.e.\ homeomorphism from the
space to itself) $\rho$ of the Euclidean $n$-sphere
$S^n$. It is known that whenever $\rho$ is conjugate, in the group of
automorphisms of $S^n$, to an orthogonal transformation, then the
fixed point set is homeomorphic to the $r$-sphere, for some $-1 \leq r
\leq n$, where $S^{-1}$ should be interpreted as the empty set. In 
general, however, the fixed points of $\rho$ need not form a sphere,
see \cite[Section I.5]{bredon} and the references cited
there. That the situation 
cannot be completely arbitrary, though, is shown by the following
result, which is one version of a classical theorem of Smith
\cite{smith}. This formulation of Smith's theorem follows e.g.\ from
\cite[Theorem III.5.1]{bredon} by passing to the second
barycentric subdivision of $\Delta$.

\begin{theorem}[Smith]\label{th:smith}
Let $\Delta$ be a finite simplicial complex which is a homology $n$-sphere
over $\Z_2$. Suppose $\Z_2$ acts simplicially on $\Delta$ in such a
way that every fixed simplex is fixed pointwise. Then, the
subcomplex induced by the fixed vertices of $\Delta$ is a homology
$r$-sphere over $\Z_2$, for some $-1 \leq r \leq n$.
\end{theorem}

\begin{remark}
More generally, the result holds if $\Z_2$ is replaced by $\Z_p$, $p$ prime,
throughout. The fact that all pseudomanifolds are orientable over
$\Z_2$, but not over $\Z_p$ in general, is the reason why $\Z_2$ plays
a prominent role in this paper, whereas $\Z_p$ does not.
\end{remark}

\begin{lemma} \label{le:pseudo}
Suppose $P$ is a finite, graded and bounded poset in which every
interval is a homology sphere over $\Z_2$. Then $P$ is a pseudomanifold. 
\begin{proof}
The diamond property is immediate, since the diamond-shaped poset is
the only graded homology sphere of length $2$. Thus, $P$ is thin. It
remains to show strong connectivity.

We argue by contradiction, so suppose that $P$ is a minimal
counterexample. The maximal chains of $\oP$ can be partitioned into
strongly connected components. By
minimality of $P$, different components have empty
intersection. Hence, $\Delta(\oP)$ is a disjoint union of at least two
pseudomanifolds. Since all pseudomanifolds have nonzero
$\Z_2$-homology in top dimension, $P$ cannot be a homology sphere over
$\Z_2$, and we have a contradiction.
\end{proof}
\end{lemma}

Now, we are in position to state and prove our main technical tool.

\begin{theorem} \label{th:main}
  Let $P$ be a poset which is Gorenstein$^\ast$ over $\Z_2$. Suppose
  that we have an 
  involutive automorphism $\nu$ of $P$. Then, the subposet of $P$
  induced by the fixed points of $\nu$ is Gorenstein$^\ast$ over $\Z_2$.
\begin{proof}
Let $F\subset P$ be the set of fixed points with the induced
order. Clearly, $\Hn$ and $\He$ are fixed by $\nu$, so $F \neq
\emptyset$. We must show that $[u,v]$ is graded and a homology sphere of top
dimension over $\Z_2$ for any interval $[u,v] \subseteq F$. 

By Theorem \ref{th:smith}, every interval in $F$ is a homology sphere
over $\Z_2$. We must still show, however, that it is graded, and that
the non-zero reduced homology group is in fact the top-dimensional one.

First, we show that every interval in $F$ is graded. Suppose, in order
to get a contradiction, that $I \subseteq F$ is a
minimal non-graded interval. By minimality of $I$, maximal
chains in $\oI$ of different lengths have empty
intersection. Thus, $\oI$ is a disjoint union of graded
posets. By Lemma \ref{le:pseudo}, all connected components of $\Delta(\oI)$ are
pseudomanifolds, and since $P$ is not graded, there are at least two
of them. Just as in the proof of Lemma \ref{le:pseudo}, this
contradicts $\Delta(\oI)$ being a $\Z_2$ homology sphere. Thus, every
interval is graded.  

Again, by Lemma \ref{le:pseudo}, every interval $[u,v] \subseteq F$ is a
pseudomanifold. Thus, its unique non-zero reduced homology group over
$\Z_2$ must be of top dimension.
\end{proof}
\end{theorem}

\section{The Bruhat order on twisted involutions}\label{se:twist}
Recall that $(W,S)$ is a Coxeter system. Suppose we have an involutive
group automorphism $\theta:W \to W$ 
which preserves $S$ as a set. In particular, $\theta$ must be a poset
automorphism of $\Br(W)$, and therefore, by Theorem \ref{th:hombergh},
be induced by an involutive automorphism of the Coxeter graph of
$W$. (As is readily checked, this indeed holds also for dihedral groups.)
\begin{definition}
The set $\twist(\theta)$ of {\em twisted involutions} with respect to
$\theta$ is defined by $\twist(\theta) = \{w \in W \mid \theta(w) = w^{-1}\}$. 
\end{definition}
We denote by $\Br(\twist(\theta))$ the subposet of $\Br(W)$ induced by
$\twist(\theta)$. When $W$ is a Weyl group, this poset plays a prominent role
in the study of related symmetric varieties, see Richardson and
Springer \cite{RS, RS2}. The said authors showed that $\Br(\twist(\theta))$
enjoys many of the nice properties associated with ordinary Bruhat
orders. In particular they proved that, in Weyl groups, $\Br(\twist(\theta))$
is graded with a certain geometrically defined rank function.

The special case $\theta = \id$ is particularly interesting. Note that
$\twist(\theta)$ is the set of involutions in this situation. We use the notation $\In(W) = \Br(\twist(\id))$. Incitti \cite{incitti, incitti2,
  incitti3} used (signed) permutation group interpretations to show
that when $W$ is of type $A$, $B$ or $D$, $\In(W)$ is
EL-shellable (hence Cohen-Macaulay) and Eulerian with rank function
being the average of 
the length and the absolute length. He conjectured that the same holds
for every Coxeter group (if $W$ is infinite, the properties should
hold for every interval in $\In(W)$). 

In the Weyl group
case, the aforementioned rank function studied by Richardson and
Springer \cite{RS} is equivalent to the one predicted by Incitti via 
a result of Carter \cite[Lemma 2]{carter}. As was pointed out in
\cite{dyer}, this equivalence does not extend to general Coxeter groups.

Below, we prove part of Incitti's conjecture for arbitrary Coxeter
groups, namely the Gorenstein$^\ast$ property over $\Z_2$ and the assertion
about the rank function. In fact, we prove similar properties for
arbitrary $\theta$. To see what remains unproved of the conjecture, recall that
if a poset is EL-shellable and Eulerian, then every interval is
homeomorphic to a top-dimensional sphere; in particular, the poset is
Gorenstein$^\ast$ over $\Z$, which is stronger than being
Gorenstein$^\ast$ over $\Z_2$.

\begin{theorem} \label{th:gorenstein}
Every interval in $\Br(\twist(\theta))$ is Gorenstein$^\ast$ over
$\Z_2$.
\begin{proof}
Choose arbitrary twisted involutions $v < w \in
\twist(\theta)$. Let $\inv:W \to W$ be the inversion map $u \mapsto
u^{-1}$. The composite map $\inv \circ \theta$ is an involutive poset
automorphism of $\Br(W)$, since $\inv$ and $\theta$ commute. Note that
its set of fixed points is 
$\twist(\theta)$. Applying Theorems \ref{th:BW} and \ref{th:main} to
$[v,w] \subseteq \Br(W)$ yields the result. 
\end{proof}
\end{theorem}

Although its existence is ensured by Theorem \ref{th:gorenstein}, it
requires some effort to actually describe the rank function of
$\Br(\twist(\theta))$. We need some notation.

\begin{definition}
The set $\iota(\theta)$ of {\em twisted identities} of $W$ with respect
  to $\theta$ is defined by $\iota(\theta) = \{w\theta(w^{-1})\mid w \in W\}$.
\end{definition}

Note that, in particular, $\iota(\id) = \{ e \}$, where $e$ is the
identity element in $W$. 

The following simple observation will prove useful later.

\begin{lemma} \label{le:identity}
If $s_1\dots s_k \in \iota(\theta)$, then $s_2\dots s_k\theta(s_1) \in
\iota(\theta)$, too. 
\begin{proof}
If $s_1\dots s_k = w\theta(w^{-1})$, then $s_2\dots s_k\theta(s_1) =
s_1w\theta((s_1w)^{-1})$. 
\end{proof}
\end{lemma}

\begin{definition}
Given $w \in W$, the {\em twisted absolute length} of $w$ with respect to
$\theta$ is denoted by $\ell^\theta(w)$ and defined as follows. Let
$s_1\dots s_k$ be any reduced expression for $w$. Then
$l = \ell^\theta(w)$ is the smallest natural number such that for some
choice of $i_1, \dots, i_l \in [k]$, we obtain $s_1\dots
\widehat{s_{i_1}} \dots 
\widehat{s_{i_l}} \dots s_k \in \iota(\theta)$. In other words,
$\ell^\theta(w)$ is the smallest number of 
elements that must be deleted from any reduced expression for $w$ in
order to obtain a twisted identity.   
\end{definition}

Since $e \in \iota(\theta)$ regardless of $\theta$, we can always
obtain a twisted identity by deleting generators in an expression. It
is not self-evident, however, that the above definition is 
independent of the choice of reduced expression for $w$. We now show
that it is.

\begin{lemma}
The twisted absolute length is well-defined.
\begin{proof}
Pick $w \in W$. It is well-known that any pair of reduced expressions
for $w$ is connected by a sequence of {\em braid moves}, each
replacing a factor $s_is_js_i \dots$ by the factor
$s_js_is_j \dots$, the length of each factor being $m(s_i,s_j)$, the order of
$s_is_j$. These factors may be interpreted as the two different
reduced expressions for the longest element, call it $y$, in
the dihedral parabolic subgroup $\langle s_i,s_j\rangle$. Thus, it
suffices to show that if $x \in \langle s_i,s_j\rangle \setminus \{y\}$ can be
obtained from one of the reduced expressions for $y$ by deleting $l$
generators, then the same holds for the other reduced expression. Now
we need only note that, in order to obtain $x$ from an arbitrary
reduced expression for $y$, it is necessary and sufficient
to delete one generator if $\ell(x)$ and
$\ell(y)$ have different parity, and two otherwise. 
\end{proof}
\end{lemma}

The following lemma seems very natural. To prove it, however, we have
to delve into some subtle properties of $\iota(\theta)$. To enhance
readability, we postpone the proof to the end of this section. 

\begin{lemma} \label{le:length}
Suppose $s_1\dots s_{k-1}\theta(s_1)$ is a reduced expression for $w
\in \twist(\theta)$. Then $\ell^\theta(s_2\dots s_{k-1}) = \ell^\theta(w)$.
\end{lemma}

Dyer \cite{dyer} showed that the absolute length $\ell^\prime(w)$ of
an element $w \in W$ is equal to the smallest number of generators
that need to be deleted in any reduced expression for $w$ in order to
obtain the identity element $e$. In other words, $\ell^\id =
\ell^\prime$. Thus, putting $\theta = \id$ in the following theorem
shows that the rank function of $\In(W)$ is the average of the
length and the absolute length, as conjectured by Incitti. In the Weyl
group case (for arbitrary $\theta$), we obtain an alternative
interpretation of the rank function defined in \cite{RS}.

\begin{theorem}
The rank of $w \in \Br(\twist(\theta))$ is $(\ell(w) + \ell^\theta(w))/2$.
\begin{proof}
Let $w$ be any twisted involution different from $e$. We already know that
$\Br(\twist(\theta))$ is graded. Thus, it suffices to show that $w$ covers
some element $v$ in $\Br(\twist(\theta))$ and either {\rm (i)}
$\ell(v) = \ell(w)-2$ and $\ell^\theta(v) = \ell^\theta(w)$, or {\rm
  (ii)} $\ell(v) = \ell(w)-1$ and $\ell^\theta(v) =
\ell^\theta(w)-1$. There are two cases:

\noindent {\em Case 1. There exists a reduced expression $s_1s_2\dots
  s_{k-1}\theta(s_1)$  for $w$}:

Let $v = s_2 \dots s_{k-1}$. Observe that $v\theta(v) =
s_1w\theta(s_1)\theta(s_1w\theta(s_1)) = s_1w\theta(w)s_1 =
s_1ww^{-1}s_1 = e$. Hence, $v$ is a twisted involution. Furthermore,
$s_1v\theta(s_1v) = w\theta(v) = wv^{-1} \neq e$, so that $s_1v \not
\in \twist(\theta)$. Similarly, $v\theta(s_1) \not \in
\twist(\theta)$, implying that $w$ covers $v$. 

Clearly, $\ell(v) = \ell(w)-2$. Lemma \ref{le:length} shows that
$\ell^\theta(v) = \ell^\theta(w)$, as desired.

\noindent {\em Case 2. No reduced expression $s_1\dots s_k$ for $w$
  satisfies $\theta(s_1) = s_k$}:

Choose a reduced expression $s_1\dots s_k$ for $w$. Suppose that
$\ell^\theta(w) = l$, and pick appropriate $i_1,\dots ,i_l \in [k]$ so
that $s_1\dots \widehat{s_{i_1}} \dots \widehat{s_{i_l}} \dots s_k  \in
\iota(\theta)$. Since $w\theta(w) = e$, we must have
$\ell(w\theta(s_1)) < \ell(w)$. Therefore, by the Exchange Property
and the fact that we are in Case 2, $w = w\theta(s_1)^2 = s_2\dots s_k
\theta(s_1)$. Repeating this argument, we see that $w = s_{i_1} \dots
s_k\theta(s_1)\dots \theta(s_{i_1-1})$. Applying Lemma
\ref{le:identity}, we may thus assume without loss of generality that
$i_1 = 1$.  

Now, let $v = s_2 \dots s_k$. Clearly, $\ell(v) = \ell(w) - 1$, and we
have just shown that $\ell^\theta(v) = \ell^\theta(w) - 1$. It remains
to prove that $v$ is a twisted involution. Since $v = w\theta(s_1)$, we
have $\theta(v)v =
\theta(s_1w)w\theta(s_1) = \theta(s_1)\theta(w)w\theta(s_1) =
\theta(s_1)^2 = e$. Thus, $v \in \twist(\theta)$, and we are done.
\end{proof}
\end{theorem} 

\subsection{Proof of Lemma \ref{le:length}}
For $w \in W$, let $J(w) = \{s_1, \dots, s_k\} \subseteq S$, where
$s_1\dots s_k$ 
is any reduced expression for $w$. The well-known fact that any two
reduced expressions for an element contain the same set of Coxeter generators
shows that $J(w)$ is unambiguously defined.

Given $J \subseteq S$, denote by $W_J = \langle J \rangle$ the
parabolic subgroup generated by $J$. It is well-known that every
(right) coset $W_Jw$ has a unique member $^Jw$ of minimal length. It
is characterized by the property that none of its reduced expressions
begins with a letter from $J$.

In order to prove Lemma \ref{le:length}, we need the following bit of knowledge
about the structure of $\iota(\theta)$:

\begin{lemma}\label{le:reduced}
If $w \in \iota(\theta)$, then there exists $x\in W$ such that $w =
x\theta(x^{-1})$ and $\ell(w) = 2\ell(x)$.
\begin{proof}
The assertion is trivial if $w = e$, and we proceed by induction over
$\ell(w)$. If there is a reduced expression of the form $s_1\dots
s_{k-1}\theta(s_1)$ for $w$, then $s_1w\theta(s_1)$ is a twisted
identity of smaller length, and we are done by induction.

Suppose that there is no such expression, i.e.\ that $s_k \neq
\theta(s_1)$ for every reduced expression $s_1\dots s_k$ for
$w$. Choose such an expression. Note that $\theta(w) = w^{-1}$; in
particular $\ell(w\theta(s_1)) < \ell(w)$. Since $w = w\theta(s_1)^2$, the
Exchange Property therefore implies $w = s_2\dots s_k \theta(s_1)$. Repeating
this argument, we find that $\ell(ws_i) < \ell(w)$ for all $i \in
[k]$, implying that $w = w_0(J(w))$, the longest element in the parabolic
subgroup $W_{J(w)}$. We will complete the proof by showing that no
twisted identity has these properties. Assume that $w =
x\theta(x^{-1})$. Let $J = J(w)$, and write $x = x_J\cdot {^J}x$ for
$x_J \in W_J$. The fact that $sw\theta(s) = w$ for all $s\in J$
implies $x_J^{-1} w \theta(x_J) = w$, so that we may assume
$x = {^J}x$. Hence, $w\theta({^J}x) = {^J}x$, implying that
$\theta({^J}x) = {^J}x$, 
since both elements must coincide with the minimal element in the
coset $W_J{^J}x$. This, however, means that $w = e$, a contradiction.   
\end{proof}
\end{lemma}

Thus, every twisted identity has a reduced expression of the form
$s_1\dots s_k \theta(s_k) \dots \theta(s_1)$. With this information,
we are ready to give the postponed proof.

\begin{proof}[Proof of Lemma \ref{le:length}]
Let $w = s_1 \dots s_{k-1}\theta(s_1)$ be as in the statement of the
lemma, and let $v = s_1w\theta(s_1) = s_2 \dots s_{k-1}$. 

Since $s_1x\theta(s_1)$ is a twisted identity whenever $x$ is, we
immediately obtain $\ell^\theta(w) \leq \ell^\theta(v)$. To prove the
other direction, we must show that if it is possible to omit $l$ generators
in the above expression for $w$ in order to yield a twisted identity,
then at most $l$ need to be deleted in the expression for $v$. This is
immediate if none or both of the initial $s_1$ and the terminal $\theta(s_1)$
are omitted. We may therefore suppose that exactly one of them is
deleted; without loss of generality, assume it to be the initial
one. In other words, we assume that $u = s_2\dots \widehat{s_{i_2}} \dots
\widehat{s_{i_l}} \dots s_{k-1}\theta(s_1) \in \iota(\theta)$. Thus, we can
obtain $u\theta(s_1)$ from $s_2\dots s_{k-1}$ by deleting $l-1$ generators.

If $\ell(u\theta(s_1)) > \ell(u)$, then the Exchange Property implies that
$u$ can be reached from our expression for $v$ by deleting $l$
generators, and we are done. Suppose now that $\ell(u\theta(s_1)) <
\ell(u)$. Applying Lemma \ref{le:reduced}, we may choose a reduced
expression $s_1^\prime\dots s_m^\prime \theta(s_m^\prime)\dots
\theta(s_1^\prime)$ for $u$. Omitting one of these generators yields
$u\theta(s_1)$. Thus, it is possible to choose some $t \in T$ such
that $u\theta(s_1)t = s_1^\prime\dots \widehat{s_i^\prime}\dots
s_m^\prime\theta(s_m^\prime) \dots \widehat{\theta(s_i^\prime)}\dots
\theta(s_1) \in \iota(\theta)$ for some $i \in [m]$. Noting that
$\ell(u\theta(s_1)t)<\ell(u\theta(s_1))$, we may invoke the Strong Exchange
Property to conclude that $u\theta(s_1)t$ can be obtained from
$s_2\dots s_{k-1}$ by deleting $l$ generators.
\end{proof}

\section{Involutions induced by graph automorphisms}\label{se:graphauto}
The topic of the previous section was fixed points of compositions of
the inversion map with group automorphisms induced by Coxeter graph
automorphisms. Theorem \ref{th:hombergh} shows that all other
automorphisms of irreducible Bruhat orders are induced by Coxeter graph
automorphisms (if $|S|\geq 3$). In this section we will study
involutive maps of the 
latter type. This class includes, in particular, all automorphisms of
Coxeter graphs of finite irreducible groups, with the exception of $D_4$.

Let $\varphi:W \to W$ be a group automorphism induced by an automorphism of
the Coxeter graph of $W$, such that $\varphi^2 = \id$. Mapping $S$ to itself,
$\varphi$ is also a poset automorphism of $\Br(W)$. Applying Theorem
\ref{th:main} we may conclude that every interval in the subposet of
fixed points is 
Gorenstein$^\ast$ over $\Z_2$. However, a stronger statement will be
proved in Theorem \ref{th:graphauto} below.

We need some preliminaries. Suppose $G$ is any group of (labelled) graph
automorphisms of the Coxeter graph of $W$ (at this stage, we do not
require $G$ to consist of involutions). Recall that for $J
\subseteq S$, $W_J$ is the parabolic subgroup generated by $J$. If
$J$ is finite, we again denote the longest element in $W_J$ by $w_0(J)$. Define
a set of symbols 
\[
\tS = \{\ts_J \mid J \subseteq S \text{ is a $G$-orbit, and $W_J$ is finite}\}.
\]
Steinberg proved the following theorem for finite Coxeter groups. The
other citations contain the general case.
\begin{theorem}[H\'ee \cite{hee}, M\"uhlherr \cite{muhlherr},
    Steinberg \cite{steinberg}]\label{th:steinberg} 
With suitably defined Coxeter relations, $\tS$ generates a Coxeter
system $(\tW,\tS)$ such that $\ts_J \mapsto w_0(J)$ defines an injective
group homomorphism $\phi:\tW \to W$ whose image is the subgroup $W^G$ of
fixed elements under the $G$-action.
\end{theorem}

\begin{remark}\label{re:cases}
The group $\tW$ in Theorem \ref{th:steinberg} can be recognized by a simple
inspection of the Coxeter graph of $W$, see \cite{crisp,crisperr}. We
do not review this procedure here. However, the following three cases
will be of particular interest to us later. They can easily be checked
by direct computation. In all three cases, the group 
acting is $\Z_2$, and it acts in the only possible, non-trivial
way. We obtain: $\widetilde{A_n} \cong B_{\lceil 
  \frac{n}{2}\rceil}$, $\widetilde{D_n} \cong B_{n-1}$ and
$\widetilde{E_6} \cong F_4$.   
\end{remark}

The next lemma is a reformulation of a lemma of Crisp
\cite{crisp}. He used it to recover Theorem \ref{th:steinberg} from
his more general results. Let $\tS^\ast$ and $S^\ast$ denote the free
monoids on the alphabets $\tS$ and $S$, respectively. 

\begin{lemma}[see Lemma 15 in \cite{crisp}] \label{le:crisp}
For $w \in W$, let $w^\ast \in S^\ast$ be a fixed reduced expression
  for $w$ (chosen arbitrarily). Then, the map $\phi^\ast:\tS^\ast
  \to S^\ast$ defined by $\ts_J \mapsto w_0(J)^\ast$ maps expressions
  that are reduced in $\tW$ to expressions that are reduced in $W$.
\end{lemma}

Since $\tW$ is a Coxeter group, one can define the Bruhat order
$\Br(\tW)$. Applying $\phi$, this gives a partial ordering on the fixed
points of $G$. It is not clear, though, whether it coincides with the
induced subposet of $\Br(W)$.  

The situation for the weak order is simple.
\begin{proposition}\label{pr:weak}
Let $F(W)$ be the subposet of the weak order on $W$ induced by the fixed point
subgroup $W^G \cong \tW$. Then, $F(W)$ is isomorphic to the weak order
on $\tW$.
\begin{proof}
In this proof, for brevity, let $P$ be the weak order on $\tW$. The
map $\phi$ defined in Theorem \ref{th:steinberg} is a bijection of
sets $P \to F(W)$. By Lemma \ref{le:crisp}, it is order-preserving. 

Consider an arbitrary ordered pair $u \leq v = uw$ in $F(W)$, where
$\ell(v) = \ell(u) + \ell(w)$. Note that $w \in F(W)$. Choose
reduced expressions $r_1$ and 
$r_2$ for $\phi^{-1}(u)$ and $\phi^{-1}(w)$, respectively. Note
that $r_1r_2$ is an expression for $\phi^{-1}(v)$ which is reduced,
too. (Otherwise, a subexpression of it would, by Lemma \ref{le:crisp}
and the Deletion Property,
be mapped by $\phi^\ast$ to an expression for $v$ shorter than
$\ell(v)$, a contradiction.) Thus, $\phi^{-1}(u) \leq \phi^{-1}(v)$ in
$P$, and we conclude that $\phi$ is a poset isomorphism.
\end{proof}
\end{proposition}

Aided by Theorem \ref{th:main}, we are able to prove the analogous result for Bruhat order when $G = 
\Z_2$. In particular, this is the only possibility if $W$ is irreducible and finite,
unless $W = D_4$. (It is easy to check that the corresponding
statement holds also for the three-element symmetry group associated
with $D_4$.) However, the result is true for any $G$; a proof was suggested to us by an anonymous referee. It is stated in Appendix \ref{se:referee}.    

\begin{theorem}\label{th:graphauto}
Let $\varphi$ be an involutive group automorphism of $W$ which
preserves $S$. Then, the subposet of $\Br(W)$ induced by the fixed
point group $W^{\{\id,\varphi \}} \cong \tW$ is isomorphic to $\Br(\tW)$.  
\begin{proof}
Denote by $F(W)$ the subposet of $\Br(W)$ induced by the fixed points
of $\varphi$. By Lemma \ref{le:crisp}, the bijection $\phi:\Br(\tW)
\to F(W)$ is order-preserving. 

Choose $w \in F(W)$. Define
$\widetilde{I} = [e,\phi^{-1}(w)] \subseteq \Br(\tW)$ and $I = [e,w]
\subseteq F(W)$. The restriction of $\phi$ to $\widetilde{I}$ is an
order-preserving injection $\widetilde{I} \to I$. Thus, on the
order-complex level, $\widetilde{I}$ is isomorphic to a subcomplex of
$I$. To show that $\phi$ is a poset isomorphism, it
suffices to show that $I$ and $\widetilde{I}$ are isomorphic as
simplicial complexes. Theorems \ref{th:BW} and \ref{th:main} show that
both complexes are pseudomanifolds. Thus, 
we are done once we have shown that the length of $\widetilde{I}$ is equal
to the length of $I$, since a pseudomanifold obviously cannot be a
proper subcomplex of another pseudomanifold of the same dimension.

Consider a saturated chain $e = v_0 < v_1 < \dots < v_k =
\phi^{-1}(w)$ in the {\em weak} order on $\tW$. By Proposition
\ref{pr:weak}, $e = \phi(v_0) < \dots < \phi(v_k) = w$ is a saturated
chain in the subposet of the weak order on $W$ induced by the fixed
points of $\varphi$. In particular, it is a chain in $I$, and it remains to
show that it is saturated. Suppose not; then we have $\phi(v_i) < x <
\phi(v_{i+1})$ for some fixed point $x$ and some $i$. By the nature of
weak order and the map $\phi$, $\phi(v_{i+1}) = \phi(v_i)w_0(J)$, for
some $J = \{s, \varphi(s)\} \subseteq S$, and $\ell(\phi(v_{i+1})) =
\ell(\phi(v_i)) + 
\ell(w_0(J))$. This implies $x = \phi(v_i)y$ for some $y \in W_J
\setminus \{e,w_0(J)\}$. By Theorem \ref{th:steinberg}, $y$ is not a
fixed point of $\varphi$, contradicting the fact that $x$ and
$\phi(v_i)$ are, and we are done.
\end{proof}
\end{theorem}

\subsection{Elements that commute with the top element}
If $W$ is a finite Coxeter group, we know that $\Br(W)$ has a top
element $w_0$. It 
is well-known (see e.g.\ \cite{BB}) that the map $W \to W$
defined by $x \mapsto w_0xw_0$ is
an automorphism of Bruhat order. Being the unique element of maximal
length, $w_0$ is clearly an involution. Hence, the above map is an involutive
automorphism. Its fixed points are the elements that commute with
$w_0$.

If $W$ is a finite, irreducible Coxeter group, it follows from the
classification of such groups that there exists a unique finite and
irreducible Coxeter group $\lW$ whose set of exponents is the set of
odd exponents of $W$. 

\begin{theorem}
Suppose $W$ is a finite, irreducible Coxeter group. The induced Bruhat
order on the set of $w_0$-commuting elements in $W$ is then isomorphic to $\Br
(\lW)$. Similarly, the induced weak order on these elements is
isomorphic to the weak order on $\lW$.
\begin{proof}
Suppose $W$ is irreducible and finite. Denote the mapping $x \mapsto
w_0xw_0$ by $\varphi$. It is well-known (see
\cite[Exercise 4.10]{BB}) that $\varphi$ is the identity
mapping iff all exponents of $W$ are odd, in which case the theorem is
trivially true.  

The dihedral case $W = I_2(m)$ is easily verified: if $m$ is even,
$\varphi$ is the identity map, and if $m$ is odd, $\varphi$ only fixes
$e$ and $w_0$. 

Now suppose $|S| \geq 3$. Since $\varphi$ is not only an automorphism
  of $\Br(W)$, but also a group automorphism of $W$, it follows from
  Theorem \ref{th:hombergh} that $\varphi$ is induced by a graph
  automorphism of the Coxeter graph of $W$. If $W$ has an even
  exponent, $\varphi$ thus coincides with the automorphism induced by the
  unique non-trivial Coxeter graph automorphism, implying that the
  fixed subgroup is isomorphic to $\tW$. The groups with an even
  exponent are $A_n$, $D_{2n+1}$ and $E_6$, and if $W$ is one of these
  groups, we have $\tW = \lW$ (see Remark \ref{re:cases}). Applying Theorem
  \ref{th:graphauto} and Proposition \ref{pr:weak} yields the claimed results.
\end{proof}
\end{theorem}
\appendix
\section{Generalizing Theorem \ref{th:graphauto}}\label{se:referee}
We are most grateful to an anonymous referee for pointing out a way to generalize Theorem \ref{th:graphauto} to arbitrary automorphism groups. In this appendix, we state the referee's argument, thereby proving the following theorem. We mantain the notation of the previous section.
\begin{theorem}\label{th:general}
Let $G$ be a group of automorphisms of $W$ that preserve $S$. Then, the subposet of $\Br(W)$ induced by the fixed point group $W^G \cong \tW$ is isomorphic to $\Br(\tW)$.   
\end{theorem}
In the proof, we will use the following characterization of the Bruhat order:
\begin{lemma}[Deodhar \cite{deodhar}, Theorem 1.1]\label{le:deodhar}
The Bruhat order is the unique partial order $\leq$ on $W$ which obeys the following two properties:
\begin{enumerate}
\item $e \leq w$ for all $w\in W$
\item given $s\in S$ and $w_1,w_2 \in W$ such that $\ell(sw_1)\leq \ell(w_1)$ and $\ell(sw_2)\leq \ell(w_2)$, we have $w_1\leq w_2 \Longleftrightarrow sw_1\leq w_2 \Longleftrightarrow sw_1\leq sw_2$. 
\end{enumerate}
\end{lemma}

\begin{proof}[Proof of Theorem \ref{th:general}]
Recall the homomorphism $\phi$ from Theorem \ref{th:steinberg}. It follows e.g.\ from Lemma \ref{le:crisp} that $\ell(\phi(\ts)w) = \ell(w)\pm \ell(\phi(\ts))$ for all $\ts \in \tS$, $w \in W^G$. 

Now choose $\ts \in \tS$ and $w_1, w_2 \in W^G$ with $\ell(\phi(\ts) w_1) = \ell(w_1) - \ell(\phi(\ts))$ and $\ell(\phi(\ts) w_2) = \ell(w_2) - \ell(\phi(\ts))$. We must show that $w_1\leq w_2 \Longleftrightarrow \phi(\ts)w_1\leq w_2 \Longleftrightarrow \phi(\ts)w_1\leq \phi(\ts)w_2$. The result then follows from Lemma \ref{le:deodhar}.

Trivially, $w_1\leq w_2 \Rightarrow \phi(\ts)w_1\leq w_2$, since $\phi(\ts)w_1 < w_1$. 

Assume $\phi(\ts) = s_1\dots s_l$, where $l = \ell(\phi(\ts))$. Since $\phi(\ts)$ is the top element in the parabolic subgroup $\langle s_1, \dots, s_l\rangle$, we have $\ell(s_i\phi(\ts)) = \ell(\phi(\ts)) - 1$ for all $i \in [l]$. Thus, we obtain the implications $s_1\dots s_lw_1 \leq w_2 \Rightarrow s_1\dots s_lw_1 \leq s_lw_2 \Rightarrow s_1\dots s_lw_1 \leq s_{l-1}s_lw_2 \Rightarrow \dots \Rightarrow s_1\dots s_lw_1 \leq s_1 \dots s_lw_2$ by repeatedly applying Lemma \ref{le:deodhar} in $\Br(W)$. We conclude that $\phi(\ts)w_1\leq w_2 \Rightarrow \phi(\ts)w_1\leq \phi(\ts)w_2$.

Finally, we again apply Lemma \ref{le:deodhar} repeatedly in $\Br(W)$ to prove the implications $s_1\dots s_l w_1 \leq s_1\dots s_l w_2 \Rightarrow s_2\dots s_lw_1 \leq s_2 \dots s_l w_2 \Rightarrow \dots \Rightarrow w_1 \leq w_2$. Thus, $\phi(\ts)w_1\leq \phi(\ts)w_2 \Rightarrow  w_1\leq w_2$, and we are done.
\end{proof}

\end{document}